\documentclass[11pt]{article}

\usepackage{a4wide,latexsym,graphicx,epstopdf,amsmath,varioref}
\DeclareGraphicsExtensions{.eps,.jpg}

\date{22 September 2017}

\title{A note on preconditioning weighted linear least squares,
       with consequences for weakly-constrained variational data assimilation}

\author{
   S. Gratton\thanks{Universit\'e de Toulouse, INP, IRIT, Toulouse, France.
            Email: serge.gratton@enseeiht.fr },
  ~S. G\"{u}rol\thanks{CERFACS,  Toulouse, France.
            Email: selime.gurol@cerfacs.fr},
  ~E. Simon\thanks{Universit\'e de Toulouse, INP, IRIT, Toulouse, France.
            Email: ehouarn.simon@enseeiht.fr}
  ~and Ph. L. Toint\thanks{NAXYS, University of Namur, Namur, Belgium.
            Email: philippe.toint@unamur.be}
}

\newcommand{\numsection}[1]{\section{#1}\setcounter{equation}{0}}

\newcommand{\beqn}[1]{\begin{equation}\label{#1}}
\newcommand{\eeqn}{\end{equation}}
\renewcommand{\Re}{\hbox{I\hskip -2pt R}}
\newcommand{\smallRe}{\hbox{\footnotesize I\hskip -2pt R}}
\newcommand{\tilA}{\tilde{A}}
\newcommand{\tilL}{\tilde{L}} 
\newcommand{\tilM}{\tilde{M}} 
\newcommand{\calB}{{\cal B}}
\newcommand{\calH}{{\cal H}}
\newcommand{\calM}{{\cal M}}
\newcommand{\calO}{{\cal O}}
\newcommand{\bpr}{{\bf Proof.} \hspace{1.5mm}}
\newcommand{\epr}{\hfill $\Box$ \vspace*{1em}}
\newcommand{\proof}[1]{
\begin{list}{}{
\setlength{\topsep}{0.0pt}
\setlength{\partopsep}{0.0pt}
\setlength{\leftmargin}{0.025\textwidth}
\setlength{\rightmargin}{0.5\leftmargin}
\setlength{\labelwidth}{0.5\leftmargin}
\setlength{\labelsep}{0.25\leftmargin}}
\item \bpr #1 \epr \noindent
\end{list}}
\newtheorem{theorem}{Theorem}[section]
\newtheorem{corollary}[theorem]{Corollary}
\newtheorem{lemma}[theorem]{Lemma}
\newcommand{\lthm}[2]{\vspace{\baselineskip} 
\noindent\framebox[\textwidth]{\parbox{0.95\textwidth}{
\begin{theorem} \label{#1} \rm #2 \end{theorem} } } \vspace{\baselineskip} }
\newcommand{\lcor}[2]{\vspace{\baselineskip} 
\noindent\framebox[\textwidth]{\parbox{0.95\textwidth}{
\begin{corollary} \label{#1} \rm #2 \end{corollary} } } \vspace{\baselineskip} }
\newcommand{\llem}[2]{\vspace{\baselineskip} 
\noindent\framebox[\textwidth]{\parbox{0.95\textwidth}{
\begin{lemma} \label{#1} \rm #2 \end{lemma} } } \vspace{\baselineskip} }
\newcommand{\req}[1]{(\ref{#1})}
\newcommand{\eqdef}{\stackrel{\rm def}{=}}
\newcommand{\tim}[1]{\;\; \mbox{#1} \;\;}
\newcommand{\sfrac}[2]{{\scriptstyle \frac{#1}{#2}}}
\newcommand{\half}{\sfrac{1}{2}}
\newcommand{\bigmax}{\displaystyle \max}

\begin{document}

\maketitle
\begin{abstract}
The effect of preconditioning linear weighted least-squares using an
approximation of the model matrix is analyzed, showing the interplay of the
eigenstructures of both the model and weighting matrices.  A small example is
given illustrating the resulting potential inefficiency of such
preconditioners. Consequences of these results in the context of the
weakly-constrained 4D-Var data assimilation problem are finally discussed.
\end{abstract}

\vspace*{1cm}

\numsection{Introduction}

Solving weighted linear least-squares problems, that is optimization problem
of the form 
\beqn{eq:wlsp}
\min_{x\in \smallRe^n} \frac{1}{2} \| Ax-b\|^2_{W^{-1}}
\eeqn
(with $A\in \Re^{m\times n}$, $b\in \Re^m$ and $W\in \Re^{m\times m}$
symmetric positive-definite), is an ubiquitous problem in applied 
mathematics (see \cite{Bjor96,GoluVanL96,KariKura04} for an introduction to
this domain and its vast literature), in particular when modelling in the
presence of uncertainties in the data and/or the model itself. The particular
application which motivates this paper is the solution of the
weakly-constrained 4D-Var problem in data assimilation
\cite{Trem06,Trem07,VidaPiacleDi04,Zupa97}, a mathematical formulation used,
among others, for weather forecasting
\cite{JanjLang15,LeeLee03,FishTremAuviTanPoli11,Bonaetal17} and oceanography
\cite{DiLoetal07}. 

It is well-known that the solution of \req{eq:wlsp} is given by solution of
the system of ``normal equations'' 
\beqn{eq:ls}
(A^TW^{-1}A)x=A^TW^{-1}b.
\eeqn
In many applications of interest (such as data assimilation in the earth
sciences), this system can be so large that the use of factorizations becomes
impractical, and one is then led to applying iterative methods, such as 
Krylov methods \cite{Saad96}. However, these methods typically require
preconditioning for achieving computational efficiency, often in the context
of parallel computing. Building preconditioners for general symmetric
positive-definite matrices has been widely investigated over the years (and is
out of the scope of the present study). It is however fair to say that the
choice of a good preconditioner is often far from easy and typically relies on
experience and on the details of the problem at hand.  When the matrix to
precondition is that of the system of normal equations \req{eq:ls} and $W$ is
known, one might consider that a reasonable preconditioner may be obtained by
using a suitable approximation of the matrix $A$.

The purpose of this short paper is to show why this strategy may sometimes be
ineffective.  While practitioners have been aware of the difficulty for some
time (see \cite{FishGuro13,FishGuro17,GratGuroSimoToin17b,JanjZhan17} for
example), a formal analysis, and hence a complete understanding, has been
missing so far.  A first step in this direction was made by Braess and Peisker
in \cite{BraePeis86}, where they showed (in a slighly different context) that,
if $A$ is square, symmetric and positive-definite, and if $W$ is the identity
matrix, then preconditioning $A^2$ (which corresponds to unweighted symmetric
least-squares) using the square of an approximation of $A$ as a preconditioner
might lead to a situation worse than not preconditioning at all, unless $A$
and its preconditioner commute. Our objective is to elaborate further and to
provide an analysis for the case where $A$ need not be symmetric nor
positive-definite while still requiring that $A$ be square and nonsingular. As
it turns out, this framework is general enough to cover our application in
data assimilation.

The paper is organized as follows.  Section~\ref{basic-s} proposes the main
analysis and relevant theorem, while a small illustrative numerical example is
presented in Section~\ref{example-s}.  The consequences of our analysis for
the weakly-constrained 4D-Var data assimilation are then discussed in
Section~\ref{da-s} and some conclusions finally drawn in
Section~\ref{concl-s}.

\numsection{Preconditioning weighted linear least squares}\label{basic-s}

Let the non-singular matrix $\tilA \in  \Re^{n\times n}$ be an approximation
(in a sense yet to be defined) of $A\in \Re^{n \times n}$.  Then the inverse
of the matrix
\beqn{eq:prec}
P=\tilA^TW^{-1}\tilA
\eeqn
may be used to construct a preconditioner for the system (\ref{eq:ls}),  yielding
\begin{equation}
\label{eq:slp}
%(\tilA^{-1}W\tilA^{-T})(A^TW^{-1}A)x=
P^{-1}(A^TW^{-1}A)x
=P^{-1}A^TW^{-1}b
= (\tilA^{-1}W\tilA^{-T})A^TW^{-1}b.
\end{equation}
The condition number of the preconditioned system matrix
$A_p=(\tilA^{-1}W\tilA^{-T})(A^TW^{-1}A)$ -- and thus the "quality" of the
preconditioner $P$ -- naturally depends on the approximation $\tilA$  and
the weight matrix $W$. A trivial (but useless)  choice is $\tilA=A$, resulting in the
condition number of $A_p$ being equal to 1.

We now show  that $\sigma((\tilA^{-1}W\tilA^{-T})(A^TW^{-1}A))$, the
spectrum of the preconditioned system matrix, is bounded by a function of the
error of $\tilA$ as an approximation of $A$ and the condition number of 
the matrix $W$.

\lthm{thm:evps}{
Let $(A,\tilA)\in \Re^{n\times n}\times\Re^{n\times n}$ be non singular
matrices, and let $W$ be a symmetric positive-definite matrix in $\Re^{n
  \times n}$ and let
\beqn{Ap-def}
A_p \eqdef (\tilA^{-1}W\tilA^{-T})(A^TW^{-1}A).
\eeqn
Then
\beqn{spectrum}
\sigma(A_p )
\subset \calB\Big(1,(1+\kappa_2(W))\| E\|_2+\kappa_2(W)\|E\|_2^2\Big)
\eeqn
where $E \eqdef A\tilA^{-1} - I_n $ the approximation error of $A$ by $\tilA$
,$\kappa_2(W)=\|W\|\,_2\|W^{-1}\|_2$ the condition number of $W$ in the
Euclidean norm and $\calB(a,r)$ is the closed ball of radius $r$ centered in $a$ .}

\proof{
We first note that, because $\tilA$ is non singular, the eigenvalues of
$A_p$ are identical to the eigenvalues of
\[
\begin{array}{lcl}
F&=&W\tilA^{-T}A^TW^{-1}A\tilA^{-1}\\
&=&W(I_n+E^T)W^{-1}(I_n+E)\\
&=&I_n+E+WE^TW^{-1}+WE^TW^{-1}E\\
&\eqdef&I_n+G.
\end{array}
\]
Now let $(\lambda,v)$ be an eigenpair of $F$, with $\|v\|_2=1$. By definition, we have that
 \begin{equation}
\label{eq:vp}
v^TFv=1+v^TGv=\lambda \| v\|_2^2=\lambda,
 \end{equation}
and therefore, using the Cauchy-Schwarz and triangle inequalities, that
\beqn{eq:ball}
\lambda\in\mathcal{B}(1,\|G \|_2).
\eeqn
Now, using $\|E\|=\|E^T\|$,
\begin{eqnarray}
\|G\|_2
&   =  & \|E+WE^TW^{-1}+WE^TW^{-1}E \|_2\\
& \leq & \|E\|_2 + \|WE^TW^{-1}\|_2 + \| WE^TW^{-1}E \|_2\\
& \leq & \|E\|_2 + \|W\|_2 \|E^T\|_2 \|W^{-1}\|_2 + \|W\|_2 \|E^T\|_2\|W^{-1}\|_2 \|E\|_2\\
& \leq & (1+\kappa_2(W))\|E\|_2+\kappa_2(W)\|E\|_2^2
\end{eqnarray}
Combining this inequality with \req{eq:ball} then gives \req{spectrum}.
}  %\end{proof}

\noindent
It results from this theorem that the condition number of $W$ and the
approximation error $E$ interact, and that a large condition number of $W$
then requires the error $\| E\|_2$ to be correspondingly small in order to
guarantee a small bound on the eigenvalues of the preconditioned system. Thus
the choice of the approximation of the system matrix $A$, and thus of the
preconditioner, should take the weighting matrix $W$ into
account, as to ensure that $\kappa(W)\| E\|_2=\calO(1)$.

Following \cite{BraePeis86}, we now define, $\kappa (D,C)$, the condition
number of a symmetric positive definite matrix $D$ with respect to a symmetric
positive definite matrix $C$ by
\beqn{eq:cond}
\kappa(D,C)=\min_{0<\gamma_1<\gamma_2}\frac{\gamma_2}{\gamma_1}
\tim{ subject to  }
\gamma_1x^TCx\leq x^TDx\leq\gamma_2 x^TCx \tim{for all} x \in \Re^n.
\eeqn

\noindent
The following easy property then follows, where $\lambda_{\min}(M)$
(resp. $\lambda_{\max}(M)$) denotes the smallest (resp. largest) eigenvalue of
the matrix $M$.

\lthm{prop:cond-prec}{
Let $(D,C)$ two symmetric positive-definite matrices. Then
\begin{equation}
\kappa(D,C)=\frac{\lambda_{\max}(C^{-1}D)}{\lambda_{\min}(C^{-1}D)}.
\end{equation}
}

\proof{
If $C^{1/2}$ its symmetric square root of $C$ and if $y = C^{1/2}x$, we
obtain that, for all $ 0<\gamma_1<\gamma_2 $, \req{eq:cond} is equivalent to
\[
\gamma_1 \|y\|_2^2 \leq  y^TC^{-1/2}DC^{-1/2}y
\leq \gamma_2 \|y\|_2^2 \tim{ for all } y \in \Re^n.
\]
The optimal  constant $\gamma_1$ (resp. $\gamma_2$) is equal to the smallest
(resp. largest) eigenvalues of the matrix  $C^{-1/2}DC^{-1/2}$ which is also
the smallest (resp. largest) eigenvalue of the matrix  $C^{-1}D$.
} % epr

\noindent
We now provide an upper bound of the condition number of $A^TW^{-1}A$ with
respect to $\tilA^TW^{-1}\tilA$.

\lcor{cor:cond_prec}{
Let $(A,\tilA)\in \Re^{n\times n}\times\Re^{n\times n}$ be non singular matrices,
and let $W$ be a symmetric positive-definite matrix.
Then, if $E =A\tilA^{-1} - I_n $ is the approximation error of $A$ by $\tilA$, 
and assuming that
\beqn{eq:error-assum}
\| E\|_2
< \frac{-(1+\kappa_2(W))+\sqrt{(1+\kappa_2(W))^2+4\kappa_2(W)}}
         {2\kappa_2(W)},
\eeqn
one has that
\[
\kappa(A^TW^{-1}A,\tilA^T W^{-1}\tilA)
\leq\frac{1+(1+\kappa_2(W))\| E\|_2 + \kappa_2(W)\| E\|_2^2}
         {1-(1+\kappa_2(W))\| E\|_2 - \kappa_2(W)\| E\|_2^2}.
\]
}

\proof{
  From Theorem~ \ref{prop:cond-prec}, one has that
  $\kappa(A^TW^{-1}A,\tilA^TW^{-1}\tilA)$ is the ratio between the largest
  and smallest eigenvalues of the matrix $A_p$. Furthermore, the assumption
  \req{eq:error-assum} guarantees that
  $1-(1+\kappa_2(W))\| E\|_2-\kappa_2(W)\| E\|_2^2>0$. The desired conclusion then
  follows from the observation that, because of Theorem \ref{thm:evps}, the
  eigenvalues of the matrix $A_p$ defined in \req{Ap-def} all belong to
  $\calB(1,(1+\kappa_2(W))\| E\|_2+\kappa_2(W)\| E\|_2^2)$. 
}%\end{proof}

\begin{figure}[t]
\begin{center}
\hspace{-2cm}\begin{minipage}[c]{7cm}
\begin{center}
\includegraphics[width=7cm,keepaspectratio=true]
{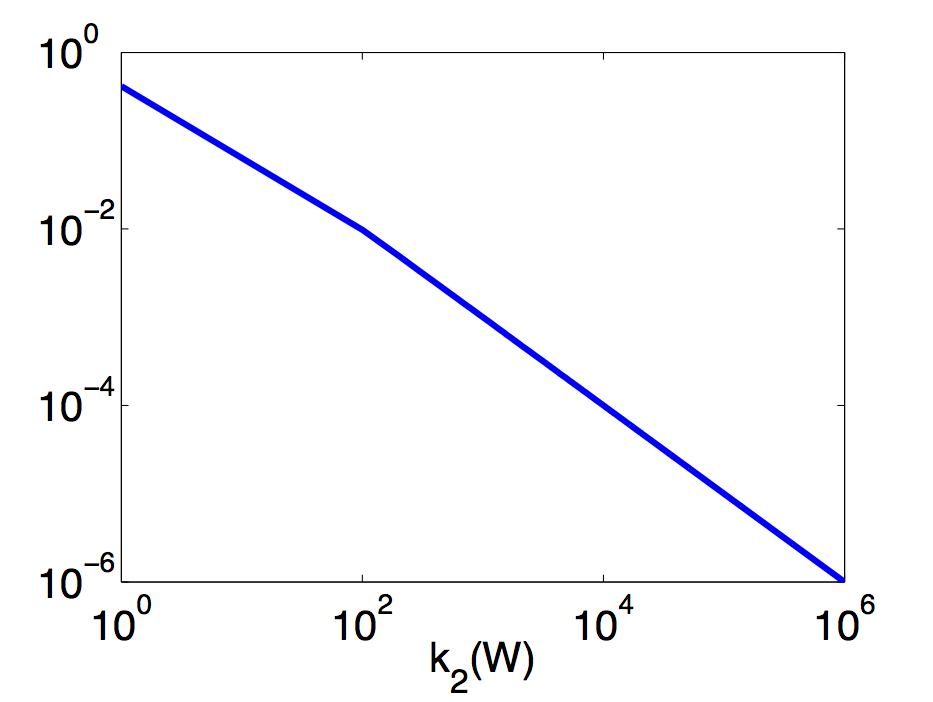}
(a) Upper bound (\ref{eq:error-assum}) on the error $\| E\|_2$ as a
    function of the condition number of $W$ (logarithmic scales). 
\end{center}
\end{minipage}
\hspace{1cm}\begin{minipage}[c]{7.0cm}
\begin{center}
\includegraphics[width=7cm,keepaspectratio=true]
{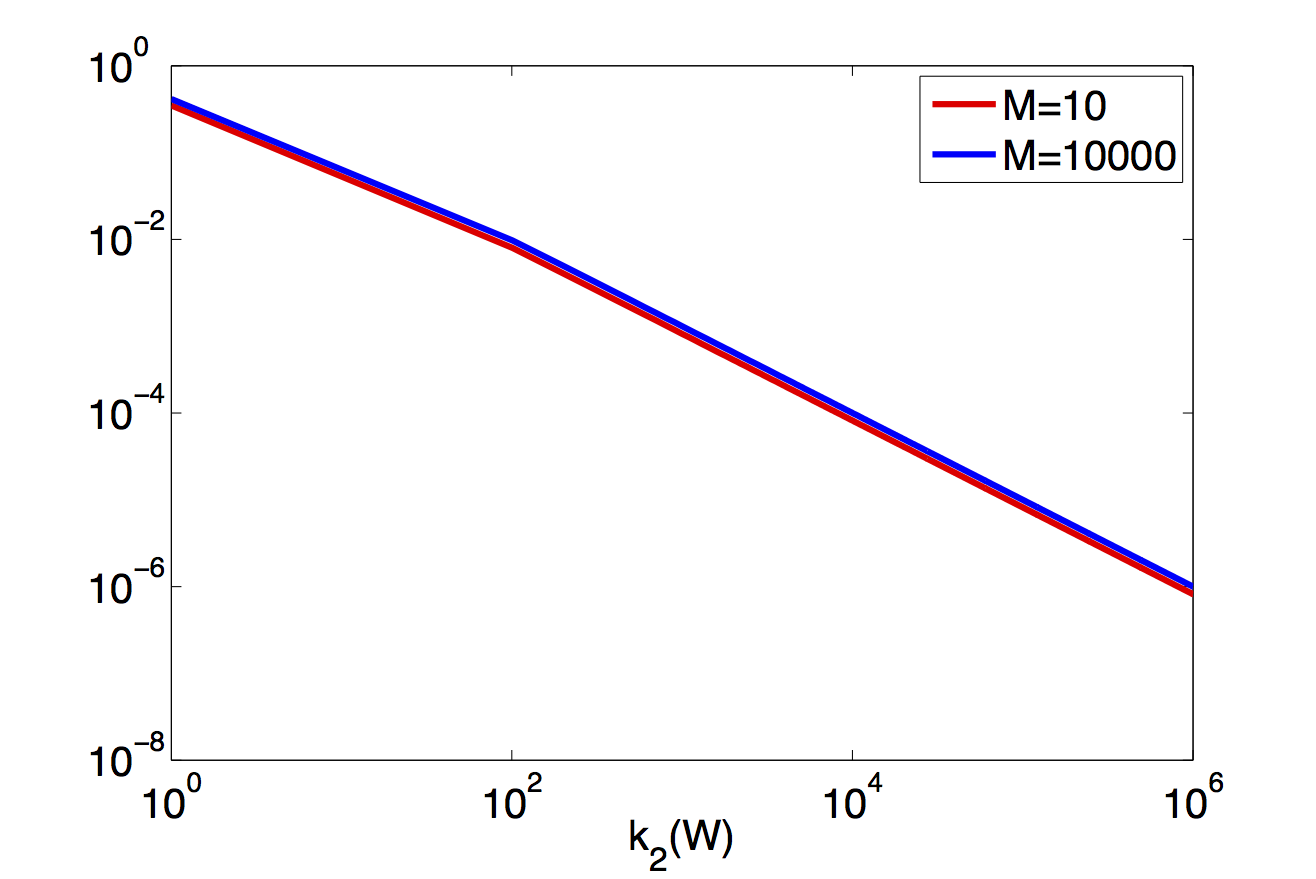}
(b) Upper bound $g$ on the error $\| E\|_2$ as a function of the
    condition number of $W$ (logarithmic scales).
\end{center}
\end{minipage}
\end{center}
\label{fig:bound_cond}  
\end{figure}

\noindent
The condition (\ref{eq:error-assum}) has a strong impact on numerical applications.
Observe that the upper bound on the error $\| E\|_2$  stated in \req{eq:error-assum}
is less than one and tends to zero when $\kappa_2(W)$ grows (see
Figure~\ref{fig:bound_cond}~(a)). For instance, a condition number $\kappa_2(W)=100$ imposes
an approximation error of the order of $10^{-2}$.  Furthermore, if one aims at a preconditioned
matrix $A_p$ with a condition number bounded above by $M>0$, then the requirement
\[
\kappa(A^TW^{-1}A,\tilA^TW^{-1}\tilA)
\leq\frac{1+(1+\kappa_2(W))\| E\|_2+\kappa_2(W)\| E\|_2^2}
         {1-(1+\kappa_2(W))\| E\|_2-\kappa_2(W)\| E\|_2^2}
\leq M
\]
results in an upper bound for the approximation error given by
\[
\| E\|_2
\leq\frac{-(1+\kappa_2(W))+\sqrt{(1+\kappa_2(W))^2+4\kappa_2(W)\frac{M-1}{M+1}}}
         {2\kappa_2(W)}
\eqdef g(\kappa_2(W),M).
\]
The evolution of $g$ with respect to $\kappa_2(W)$ is shown in
Figure~\ref{fig:bound_cond}~(b) for two values of $M$. We note that even
relatively large bounds on the condition number of $A^TW^{-1}A$ with respect
to $\tilA^TW^{-1}\tilA$ impose small approximation errors, especially when $W$
has a large condition number.

\numsection{A simple illustrative example}\label{example-s}

We now illustrate the impact of the preconditioners \req{eq:prec} on the
eigenvalues and condition number of the preconditioned system matrix
\req{eq:slp} in a very simple case. Let $\alpha\geq1$ be a parameter
corresponding to the condition number of the weight matrix $W$. We define
\[
A=\left(\begin{array}{cc}1&0 \\ \alpha&1  \end{array}\right)
\tim{ and }
W =\left(\begin{array}{cc}\alpha&0 \\ 0&1  \end{array}\right).
\]
It can then be verified that the matrices $A^TA$ and $A^TW^{-1}A$ both
have their condition numbers tending to infinity when $\alpha$ grows.
We now introduce the approximation of $A$ given by
\[
\tilA= \left(\begin{array}{cc}1&0 \\ \alpha+2&1  \end{array}\right).
\]
It is now possible to construct a "good preconditioner" $P^{-1}=\tilA^{-1}\tilA^{-T}$
of the matrix $A^TA$ in the sense that, while the condition number of $A^TA$  goes
towards infinity when $\alpha$ grows, the condition number of the matrix $A^TA$
with respect to  $\tilA^{T}\tilA$ is constant. In this specific case, one has
$\kappa(A^TA,\tilA^T\tilA)$ which is the same as $\kappa_2(\tilA^{-T}
({A}^{T}A)\tilA^{-1})$ is approximately equal to 33.9706.

However, the approximation error $E$ associated with this matrix is
\[
E=A\tilA^{-1} - I_2
= \left(\begin{array}{rr}0&0 \\ -2&0  \end{array}\right),
\]
leading to $\| E\|_2=2$. From Theorem~\ref{thm:evps}, one then has that the
eigenvalues of $A_p= (\tilA^{-1}W\tilA^{-T})(A^TW^{-1}A)$
belong to the closed ball $ \calB(1,1+6\alpha)$, which makes it possible for the largest
eigenvalue to tend to infinity with $\alpha$. Indeed, this is is what happens
in this example. One has that
\[
A_p=\left(\begin{array}{cc}1-2\alpha^2          & -2\alpha \\
                           2\alpha^3+4\alpha^2-2&  2\alpha^2+4\alpha+1
          \end{array}\right).
\]
It can be shown that the eigenvalues of $A_p$ are
$1+2\alpha\pm 2 \sqrt{\alpha(\alpha+1)}$, and so  the largest one tends to infinity
when $\alpha$ grows. Moreover, $\kappa(A^TW^{-1}A,\tilA^TW^{-1}\tilA)$, the
condition number of $A^TW^{-1}A$ with respect to
$\tilA^TW^{-1}\tilA$, 
%$$\kappa_2(\mathcal{A})=\frac{2\alpha+1+ 2\sqrt{\alpha(\alpha+1)}}{2\alpha+1- 2\sqrt{\alpha(\alpha+1)}} \rightarrow \infty$$
therefore also tends to infinity with $\alpha$.

However, if we now define the approximation of $A$ by
\[
\tilA= \left(\begin{array}{cc}1&0 \\ \alpha+\alpha^{-1}&1  \end{array}\right),
\]
the approximation error then becomes
\[
E= \left(\begin{array}{cc}0&0 \\ \alpha^{-1}&0  \end{array}\right),
\]
leading to $\|E\|_2=\alpha^{-1}$ and $\kappa_2(W)\| E\|_2=1$. Again,
Theorem~\ref{thm:evps} says that the eigenvalues of $A_p$ belongs to
$\calB(1,1+2\alpha^{-1})$, but now the radius of this ball tends to one when
$\alpha$ grows, which results in bounded eigenvalues.
This can easily be verified as, in this case,
\[
A_p=\left(\begin{array}{cc} 1-\alpha              &          -1        \\
                            \alpha^2-\alpha^{-1}+1 &\alpha+\alpha^{-1}+1
          \end{array}\right).
\]
which has two distinct eigenvalues $\half(2+\alpha^{-1}\pm\sqrt{4\alpha^{-1}+\alpha^{-2}})$
tending to one when $\alpha$ grows, as does
$\kappa(A^TW^{-1}A,\tilA^TW^{-1}\tilA)$.

\numsection{Application to weakly-constrained data assimilation}\label{da-s}

We now turn to the implications of the above results for our motivating
application, the weakly-constrained 4D variational formulation for data
assimilation. In this context, one attempts to fit an initial state $x_0$ so
as to fit observations $y_j$ taken from the evolution of a dynamical model
$\calM$ over $N_{sw}$ time windows.  We refer to \cite{Trem06,Trem07}
for further details and motivation for this formulation, but, for our
present purposes, it is enough to know that it involves the (often
approximate) solution of the optimization problem
\beqn{eq:wk4DVar}
\min_{{\bf x} \in\Re^{n}}  
     \frac{1}{2} \|x_0 - x_b\|_{\ B^{-1}}^2
   + \frac{1}{2} \sum_{j=0}^{N_{sw}}
     \left\| {\cal H}_j \big( x_j\big) - y_j \right\|_{{\bf R}_j^{-1}}^2
   +\frac{1}{2}\sum_{j=1}^{N_{sw}}
     \|x_j - \calM_j(x_{j-1}) \|_{{\bf Q}_j^{-1}}^2
\eeqn
where
\begin{itemize}
\vspace*{-1.2mm}
\item ${ x} = ( x_0, x_1, \ldots, x_{N_{sw}})^T \in\Re^n$
is the control variable (with $x_j = x(t_j)$),
\vspace*{-1.2mm}
\item  $x_b$ is the background given at the initial time ($t_0$),
\vspace*{-1.2mm}
\item  $y_j\in\Re^{m_j}$  is the observation vector over a given time interval,
\vspace*{-1.2mm}
\item $\calH_{j}$ maps the state vector ${x_j}$ from model space to observation space,
\vspace*{-1.2mm}
\item  $\calM_{j}$ represents an integration of the numerical model from time $t_{j-1}$ to 
$t_j$,
\vspace*{-1.2mm}
\item $B$, $R_j$ and $Q_j$ are the covariances  of the  background,
  observation and model error. 
\vspace*{-1.2mm}
\end{itemize}
This general unconstrained nonlinear least-squares problem is typically solved by
applying the Gauss-Newton algorithm, which iteratively proceeds by linearizing
$\calH$ and $\calM$ at the current iterate and then, again approximately, minimizing
the resulting quadratic function. If the operators $M_j$ are the linearized $\calM_j$
and $H_j$ are the linearized $\calH_j$, then the problem can be expressed in terms of
$\delta x = x-x_0$ as
\[
\min_{ \delta x \in \Re^{n}} \;
       \frac{1}{2} \| L \delta x - b \|_{D^{-1}}^2
   +   \frac{1}{2} \| H \delta x - d \|_{R^{-1}}^2
\]
where
\begin{equation}
\label{eq:Ldef}
L = \left(\begin{array}{ccccc}
               I_n  &    &  & &\\
            -M_1  &  I_n &  & & \\
                  & -M_2  & I_n & &\\
                  &       & \ddots  & \ddots &\\       
                  &       &         & -M_{N_{sw}} & I_n\\
                  \end{array}
\right)
\end{equation}
for suitable vectors
\[
d = ( d_0, d_1, \ldots, d_{N_{sw}})^T
\mbox{ and }
b = (   b, c_1, \ldots, c_{N_{sw}})^T,
\]
and where
\[
H = \mbox{diag}( H_0, H_1, \ldots, H_{N_{sw}}),
\;
D = \mbox{diag}( B, Q_1, \ldots, Q_{N_{sw}})
%\]
\tim{and}
%\[
R = \mbox{diag}( R_0, R_1, \ldots, R_{N_{sw}}).
\]
This particular form of the problem is called the "state formulation" and
its optimality conditions amount to (approximately) solving  linear systems of
the form 
\begin{equation}
\label{eq:state}
(L^TD^{-1}L+H^TR^{-1}H)\delta x = L^TD^{-1}b+H^TR^{-1}d
\end{equation}
An alternative, called the ``forcing formulation'', is also
possible by rewriting the problem in terms of $\delta p=L\delta x$, but we do
not consider it here because it is not amenable to parallel computation.
Its conditioning has been studied in \cite{ElSa15,ElSaNichLawl17}.

It is traditionally assumed that the term $L^TD^{-1}L$  (called the background
term) dominates in the system matrix, which then leads to preconditioners of
the form  
\beqn{eq:prec_state}
P^{-1}=  {\tilL}^{-1}\,D\,\,{\tilL}^{-T},
\eeqn
with $\tilL$ an approximation of the matrix $L$ (see \req{eq:Ldef}). This
approximation is often built by replacing in $L$ the operators $M_j$
associated with the numerical model by approximations $\tilM_j$. While the
matrix-vector product with $L$ can be done in parallel, the preconditioner
\req{eq:prec_state} involves $\tilL^{-1}$, whose parallelization potential
crucially depends on the choice of the operators $\tilM_j$. Two very simple
approximations are commonly chosen in practice: $ \tilM_j=0$ or $
\tilM_j=I_n$. The preconditioned system matrix is then $({\tilL}^{-1}{
  D}{\tilL}^{-T})(L^TD^{-1}L+H^TR^{-1}H)$.  In what follows, we focus on the
preconditioned background term $({\tilL}^{-1}D{\tilL}^{-T})(L^TD^{-1}L)$ and
we investigate the consequences of Theorem~\ref{thm:evps} for this matrix.

We first analyse the form of the approximation error
$E=L\tilL^{-1}-I_{nN_{sw}}$.  

\llem{cor:AN_err}{
  Let $L$, $\tilL$, $M_j$, $\tilM_j$ and $E=L\tilL^{-1}-I_{nN_{sw}}$.
  Then
\begin{itemize}
\item if $\tilM_j=0$, one has that $E= \left(\begin{array}{ccccc}
               0  &    &  & &\\
            -M_1  &  0 &  & & \\
                  & -M_2  & 0 & &\\
                  &       & \ddots  & \ddots &\\       
                  &       &         & -M_{N_{sw}} & 0\\
                  \end{array}
\right)$;
\item if $\tilM_j=I_n$, one has that $E= \left(\begin{array}{ccccc}
               0  &    &  & &\\
           I_n-M_1  &  0 &  & & \\
             I_n-M_2     & I_n-M_2  & 0 & &\\
                  &       & \ddots  & \ddots &\\       
         I_n-M_{N_{sw}}        &    I_n-M_{N_{sw}}   &\cdots         & I_n -M_{N_{sw}} & 0\\
                  \end{array}
  \right)$.               
\end{itemize}
}

\proof{
It can be verified that $E$ is block-lower
triangular with null blocks on the diagonal. Furthermore, one has that, for
all indeces $(i,j)$ such that $1\leq j<i\leq N_{sw}+1$,
\[
E_{i,j}
=\left\{\begin{array}{ll}(\tilM_{i-1}-M_{i-1})\tilM_{i-2}\cdots \tilM_{j} &\quad \mbox{ if } j<i-1, \\
\tilM_{i-1}-M_{i-1} & \quad \mbox{ if } j=i-1,
\end{array}\right.
\]
where $E_{i,j}\in \Re^{n\times n}$ is the $(i,j)$-th block of $E$. The
conclusions of the lemma then follow by specializing $M_j$.
} %epr

\noindent
Using those expressions for the approximation error (of $\tilM_j$ as an approximation of $M_j$),
we may then derive the following conclusions from Theorem~\ref{thm:evps}, in terms of
$\sigma_{\max}(M_j)$, the largest singular value of the linearized model matrix $M_j$.

\lcor{cor:AN_vp}{
Let $L$ and $M_j$ be defined in \req{eq:Ldef}, and let $\tilL$ be the
approximation of $L$ defined from $\tilM_j\in \{0,I_n\}$ for $j=1,\ldots,N_{sw}$.
Let $A_p=({\tilL}^{-1}D{\tilL}^{-T})(L^TD^{-1}L)$ be the preconditioned background matrix.
Then
\[
\sigma(A_p) \subset \calB(1,(1+\kappa_2(D))\rho+ \kappa_2(D)\rho^2)
\]
where
\[
\rho = \left\{\begin{array}{ll}
\bigmax_{j=1,\ldots,N_{sw}}\sigma_{\max}(M_j) & \tim{if } \tilM_j=0 \;\;(j=1\ldots,N_{sw}), \\*[2ex]
\sqrt{\frac{(nN_{sw}+1) (nN_{sw}+2)}{2} }\left[\bigmax_{j=1,\ldots,N_{sw}}\sigma_{\max}(I_n-M_j)\right]
&\tim{if } \tilM_j=I_n \;(j=1,\ldots,N_{sw}).
\end{array}\right.
\]
}

\proof{
\begin{enumerate}
\item Consider first the case where $M_j=0$. From Lemma~\ref{cor:AN_err}, we
  deduce that $E^TE$ is block diagonal and 
  \[
  E^TE=\mbox{diag}(M_1^TM_1, M_2^TM_2,\cdots,M_{N_{sw}}^TM_{N_{sw}},0).
  \]
  This then implies that $\| E\|_2=\max_{j=1,\ldots,N_{sw}}(\sigma_{\max}(M_j))$ and we can
  conclude by applying Theorem~\ref{thm:evps}.
\item  If $M_j=I_n$, then, from Corollary~\ref{cor:AN_err}, one has that $E=ST$ 
  with
  \[
  S=\left(\begin{array}{ccccc}
               0  &    &  & &\\
            I_n-M_1  &  0 &  & & \\
                  & I_n-M_2  & 0 & &\\
                  &       & \ddots  & \ddots &\\       
                  &       &         &I_n -M_{N_{sw}} & 0\\
                  \end{array}
    \right)
  \]
  and $T$ is the lower triangular matrix with the lower entries equal to one. 
  As in the previous case, one obtains that
  $\|S\|_2=\max_{j=1,\ldots,N_{sw}}\sigma_{\max}(I_n-M_j)$.
  Hence the desired conclusion follows from applying Theorem~\ref{thm:evps},
  and using the bound
  \[
  \| E\|_2
  \leq \| S \|_2  \| T\|_2
  \leq \|S\|_2\|T\|_F
  = \|S\|_2 \sqrt{\frac{(nN_{sw}+1) (nN_{sw}+2)}{2} }.
  \]
\end{enumerate} 
}%\end{proof}

\noindent
We immediately see that obtaining a well-conditioned matrix $A_p$ requires
specific assumptions on the dynamical models within a sub-window.  Choosing
$\tilM_j = 0$ will work well if the model itself is close to zero, which may
be unrealistic in many situations.  The choice $\tilM_j = I_n$ is often more
sensible if the dynamics of the model may remain limited, especially if the
time sub-windows are short.  This can be viewed as a motivation to choose
$N_{sw}$ large, but one nevertheless should remember that the gain in making
the singular value closer to 1 is offset by the dependence on the square root
term in part 2 of Corollary~\ref{cor:AN_vp}.  Obviously, the quality of the
preconditioner may improve with the quality of $\tilM_j$ as an approximation
of $M_j$, but it remains challenging to select good approximations which
preserve efficient parallel computation of $\tilL^{-1}$ (see
\cite{GratGuroSimoToin17} for an approach of this question). One should also
remember that our analysis merely provides bounds on the conditioning, which
are pessimistic by nature, and that the observation term $H^TR^{-1}H$ (which
we ignored here) may not always be negligible. The situation is therefore
often problem dependent, as has been demonstrated in
\cite{GratGuroSimoToin17b} where very different behaviours (good and bad) were
observed for two contrasting data assimilation problems.

\numsection{Conclusions}\label{concl-s}

We have provided a formal analysis of the preconditioning efficiency for
nonsingular weighted least-squares, thereby extending previous results by
Braess and Peisker \cite{BraePeis86} and vindicating the numerical experience
of several practitioners. We have also specialized the analysis to the state
formulation of the weakly-constrained data assimilation problems, an important
computational tool in the earth sciences.  While the conditioning bounds
discussed in this paper remain indicative as all bounds are, they nevertheless
provide some guidance on how to construct good parallelizable preconditioners,
a task which remains for now a problem-dependent exercize.

%  References

{\footnotesize

%\bibliography{/home/pht/bibs/refs}
%\bibliographystyle{plain}
}

\end{document}